\documentclass[11pt,reqno]{amsart}
\usepackage{amsmath}
\usepackage{amssymb}

\newtheorem{thm}{Theorem}[section]
\newtheorem{lem}{Lemma}[section]
\newtheorem{defn}{Definition}[section]

\numberwithin{equation}{section}
\newtheorem{rmk}{Remark}[section]

\def\pf{{\textit {Proof:} }}

\newcommand{\mysection}[1]{\section{#1}\setcounter{equation}{0}}

\newfont{\bb}{msbm10 at 12pt}

\def\R{\hbox{\bb R}}

\newcommand{\bal}{\begin{aligned}}      \newcommand{\eal}{\end{aligned}}
\newcommand{\ba}{\begin{array}}      \newcommand{\ea}{\end{array}}
\newcommand{\bc}{\begin{center}}     \newcommand{\ec}{\end{center}}
\newcommand{\be}{\begin{enumerate}}  \newcommand{\ee}{\end{enumerate}}
\newcommand{\beq}{\begin{eqnarray}}  \newcommand{\eeq}{\end{eqnarray}}
\newcommand{\beQ}{\begin{eqnarray*}} \newcommand{\eeQ}{\end{eqnarray*}}
\newcommand{\bi}{\begin{itemize}}    \newcommand{\ei}{\end{itemize}}
\newcommand{\bt}{\begin{tabular}}    \newcommand{\et}{\end{tabular}}
\newcommand{\bdm}{\begin{displaymath}} \newcommand{\edm}{\end{displaymath}}

\def\qed{\hfill{Q.E.D.}\smallskip}

\newcommand{\ls}{\setlength{\baselineskip}{12pt}
                 \setlength{\parskip}{3mm}}

\title{Positive energy theorem for (4+1)-dimensional asymptotically anti-de Sitter spacetimes}
\author{Yaohua Wang, Xu Xu}
\address[Yaohua Wang]{Institute of Mathematics, Academy of Mathematics and
System Sciences, Chinese Academy of Sciences, Beijing 100080, PR China
}
\email{wangyaohua@amss.ac.cn}
\address[Xu Xu]{School of Mathematics and Statistics, Wuhan University, Wuhan 430072, PR China
}
\email{xuxu2@whu.edu.cn}
\date{}

\begin{document}
\maketitle

\begin{abstract}
We define the total energy-momenta for $(4+1)$-dimensional asymptotically anti-de Sitter spacetimes, and prove
the positive energy theorem for such spacetimes.
\end{abstract}

\mysection{Introduction}\ls

Positive mass theorem is an important result in general relativity. When the cosmological constant is zero and
the spacetime is asymptotically flat, the definition of total
energy and total linear momentum was first given by Arnowitt, Deser and Misner \cite{ADM}.
Physicists then conjectured that the total mass of an isolated gravity system should be nonnegative. This conjecture
was first proved by Schoen and Yau \cite{SY1,SY2,SY3} and then by Witten \cite{Wi}. For the wide applications
in geometry and physics, higher dimensional positive mass theorems for asymptotically flat manifolds were
extensively studied \cite{B, Z, Ding, EHLS}.

When the cosmological constant is negative and the initial data is asymptotically hyperbolic, the corresponding
positive mass theorems were studied by many authors. When the second fundamental form is zero, such positive mass theorems
were obtained under different conditions \cite{CH, CN, Wax}.
When the second fundamental form is nonzero, such positive energy theorems were given in \cite{Ma, XZ, WXZ}.
For higher dimensional asymptotically anti-de Sitter spacetimes, Chru\'{s}ciel, Maerten and Tod \cite{CMT} gave a definition for the total
energy and other conserved quantities. Under certain assumptions on coordinate transformations, they obtained
some inequalities of the total energy. The first author, Xie and Zhang \cite{WXZ} obtained more general inequalities
of the total energy, without the assumptions on coordinate transformations.

In this paper, we establish similar inequalities involving total energy-momenta for $(4+1)$-dimensional asymptotically anti-de Sitter
spacetimes. We give the explicit form of the imaginary Killing spinors on the 0-slice of the anti-de Sitter spacetime after fixing a suitable Clifford representation.
We define the total energy and momenta for asymptotically anti-de Sitter initial data, and finally,
we provide the lower bound of the total energy in terms of the total momenta, which
implies the positive energy theorem.

This paper is organized as follows: In Section 2, we give the explicit form of imaginary Killing spinors on
the 0-slice of the anti-de Sitter spacetime.
In Section 3, we give the definition of total energy and total momenta.
In Section 4, we prove our positive energy theorem.
In Appendix, we provide the explicit form of the Killing vectors of the anti-de Sitter spacetime.

\mysection{The anti-de Sitter space-time}\ls

The anti-de Sitter ($AdS$) spacetime with negative cosmological constant $\Lambda$,
denoted by $(N,\widetilde{g}_{AdS})$, is a static spherically
symmetric solution of the vacuum Einstein equations.
The $(4+1)$-dimensional anti-de Sitter spacetime is indeed the hyperboloid
\begin{equation*}
\eta_{\alpha\beta}y^\alpha y^\beta=\frac{6}{\Lambda},\ \ \ \ \Lambda=-6\kappa^2(\kappa>0)
\end{equation*}
in $\mathbb{R}^{4,2}$ with the metric
\begin{equation}\label{metric of R4,2}
ds^2=-(dy^0)^2+\sum_{i=1}^{4}(dy^i)^2-(dy^5)^2.
\end{equation}
Under suitable choice of coordinates, the metric of $(4+1)$-dimensional anti-de Sitter spacetime can be written as
\begin{equation*}
\begin{aligned}
ds^2=-\cosh^2(\kappa r)dt^2+dr^2+\frac{\sinh^2(\kappa r)}{\kappa^2}\bigg(d\theta^2+\sin^2\theta (d\psi^2+\sin^2\psi d\varphi^2)\bigg).
\end{aligned}
\end{equation*}
The $t$-slice $(\mathbb{H}^4,\breve{g})$ is the hyperbolic
4-space with constant sectional curvature $-\kappa^2$.

Let the orthonormal frame of the AdS spacetime be
$$\breve{e}_0=\frac{1}{\cosh(\kappa r)}\frac{\partial}{\partial t},\ \
\breve{e}_1=\frac{\partial}{\partial r}, \  \ \breve{e}_2=\frac{\kappa}{\sinh(\kappa r)}\frac{\partial}{\partial \theta}, \ $$
$$
\breve{e}_3=\frac{\kappa}{\sinh(\kappa r)\sin\theta}\frac{\partial}{\partial \psi},\ \
\breve{e}_4=\frac{\kappa}{\sinh(\kappa r)\sin\theta \sin \psi}\frac{\partial}{\partial \varphi},
$$
and $\breve{e}^\alpha$ be the dual coframe of $\breve{e}_\alpha$.

Recall that the fifteen Killing vectors
\beQ
\begin{aligned}
 U_{\alpha\beta}=y_\alpha
\frac{\partial}{\partial y^\beta}-y_\beta \frac{\partial}{\partial y^\alpha}
 \end{aligned}
 \eeQ
generate rotations for $\R^{4,2}$ with the metric (\ref{metric of R4,2}).
By restricting these vectors to the hyperboloid
$\{\eta_{\alpha\beta}y^\alpha y^\beta=\frac{6}{\Lambda}\}$ with the induced metric,
the Killing vectors of $AdS$ spacetime can be derived.
See Appendix for the explicit form of $U_{\alpha\beta}$ on the 0-slice.

Let $\mathbb{S}$ be the spinor bundle of $(N,\widetilde{g}_{AdS})$ and its restriction to $\mathbb{H}^4$.
The spinor $\Phi_0 \in \Gamma(\mathbb{S})$ is called an imaginary Killing spinor along $\mathbb{H}^4$ if it satisfies
\begin{equation} \nonumber\\
\nabla^{AdS}_X \Phi_0+\frac{\kappa\sqrt{-1}}{2}X\cdot\Phi_0=0
\end{equation}
for each $X$ tangent to $\mathbb{H}^4$.

For the following application, we fix the following Clifford representation throughout this paper:
\begin{equation*}
\breve{e}_0 \mapsto \begin{pmatrix}
\ I & \  \\ \ &  -I
\end{pmatrix}, \ \
\breve{e}_1 \mapsto \begin{pmatrix}
\ &   I  \\ -I & \
\end{pmatrix},\ \
\breve{e}_2 \mapsto \begin{pmatrix}
\ &   i  \\ i  &  \
\end{pmatrix},
\end{equation*}
\begin{equation}\label{repre}
\ \breve{e}_3 \mapsto \begin{pmatrix}
\ &   j  \\ j  &  \
\end{pmatrix}, \ \
\ \breve{e}_4 \mapsto \begin{pmatrix}
\ &   k  \\ k  &  \
\end{pmatrix},
\end{equation}
where
\begin{equation*}
i =\begin{pmatrix}
\ \sqrt{-1} & \  \\ \ &  -\sqrt{-1}
\end{pmatrix}, \ \
j =\begin{pmatrix}
\ &   1 \\ -1 & \
\end{pmatrix},\ \
k =\begin{pmatrix}
\ &   \sqrt{-1}  \\ \sqrt{-1}  &  \
\end{pmatrix}.
\end{equation*}

Under this representation, we have
\begin{lem}
The imaginary Killing spinors along $\mathbb{H}^4$ are all of the form
\begin{equation}\label{ik}
\Phi_0=\begin{pmatrix}
u^+e^{\frac{\kappa r}{2}}+u^-e^{-\frac{\kappa r}{2}}\\
v^+e^{\frac{\kappa r}{2}}+v^-e^{-\frac{\kappa r}{2}}\\
\sqrt{-1}u^+e^{\frac{\kappa r}{2}}-\sqrt{-1}u^-e^{-\frac{\kappa r}{2}} \\ \sqrt{-1}v^+e^{\frac{\kappa r}{2}}-\sqrt{-1}v^-e^{-\frac{\kappa r}{2}}
\end{pmatrix},
\end{equation}
where \begin{eqnarray*}
u^+
&=&\Big(\lambda_1 e^{\frac{-\sqrt{-1}}{2}\varphi}\cos\frac{\psi}{2}+
    \lambda_2 e^{\frac{\sqrt{-1}}{2}\varphi}\sin\frac{\psi}{2}\Big)
   \cos\frac{\theta}{2}\\
&&
   +\Big(\lambda_3 e^{\frac{-\sqrt{-1}}{2}\varphi}\cos\frac{\psi}{2}+
   \lambda_4 e^{\frac{\sqrt{-1}}{2}\varphi}\sin\frac{\psi}{2}\Big)
   \sin\frac{\theta}{2},\\
u^-
&=&-\sqrt{-1}\Big(\lambda_1 e^{\frac{-\sqrt{-1}}{2}\varphi}\cos\frac{\psi}{2}+
   \lambda_2 e^{\frac{\sqrt{-1}}{2}\varphi}\sin\frac{\psi}{2}\Big)
   \sin\frac{\theta}{2}\\
&&
   +\sqrt{-1}\Big(\lambda_3 e^{\frac{-\sqrt{-1}}{2}\varphi}\cos\frac{\psi}{2}+
   \lambda_4 e^{\frac{\sqrt{-1}}{2}\varphi}\sin\frac{\psi}{2}\Big)
   \cos\frac{\theta}{2},\\
v^+
&=&\sqrt{-1}\Big(-\lambda_1 e^{\frac{-\sqrt{-1}}{2}\varphi}\sin\frac{\psi}{2}+
    \lambda_2 e^{\frac{\sqrt{-1}}{2}\varphi}\cos\frac{\psi}{2}\Big)
   \cos\frac{\theta}{2}\\
&&
   +\sqrt{-1}\Big(\lambda_3 e^{\frac{-\sqrt{-1}}{2}\varphi}\sin\frac{\psi}{2}-
   \lambda_4 e^{\frac{\sqrt{-1}}{2}\varphi}\cos\frac{\psi}{2}\Big)
   \sin\frac{\theta}{2},\\
v^-
&=&\Big(\lambda_1 e^{\frac{-\sqrt{-1}}{2}\varphi}\sin\frac{\psi}{2}-
    \lambda_2 e^{\frac{\sqrt{-1}}{2}\varphi}\cos\frac{\psi}{2}\Big)
   \sin\frac{\theta}{2}\\
&&
   +\Big(\lambda_3 e^{\frac{-\sqrt{-1}}{2}\varphi}\sin\frac{\psi}{2}-
   \lambda_4 e^{\frac{\sqrt{-1}}{2}\varphi}\cos\frac{\psi}{2}\Big)
   \cos\frac{\theta}{2}.\end{eqnarray*}
Here  $\lambda_1$, $\lambda_2$, $\lambda_3$, and $\lambda_4$ are
 arbitrary complex numbers.\end{lem}

\mysection{Definitions}\ls

Suppose that $N$ is a spacetime with the metric $\widetilde{g}$ of signature $(-1,1,1,1,1)$, satisfying the Einstein field
equations
\begin{equation}
\widetilde{Ric}-\frac{\widetilde{R}}{2}\widetilde{g}+\Lambda \widetilde{g}=T,
\end{equation}
where $\widetilde{Ric}$, $\widetilde{R}$
are the Ricci and scalar curvatures of $\widetilde{g}$ respectively, $T$ is the energy-momentum tensor of matter, and $\Lambda$ is the cosmological
constant.  For orthonormal frame $\{e_\alpha\}$ with $e_0$ timelike, the dominant energy condition
\begin{equation} \label{dec}
T_{00}\geq \sqrt{\sum_iT_{0i}^2} ,\ \ T_{00}\geq |T_{\alpha\beta}|
\end{equation}
is satisfied.
Let $M$ be a 4-dimensional
spacelike hypersurface in $N$ with the induced metric $g$ and $p$ be the second fundamental form of $M$ in $N$.

\begin{defn} An initial data set $(M, g, p)$ is asymptotically $AdS$ of order $\tau >2$ if\\
$(1)$ there is a compact set $K$ such that $M_\infty=M\setminus K$ is diffeomorphic to $\mathbb{R} ^4 -B_r$,
where $B _r$ is the closed ball of radius $r$ with center at the coordinate origin;\\
$(2)$ Under the diffeomorphism,
$g _{ij}=g(\breve{e} _i,\breve{e} _j)=\delta_{ij}+a_{ij}$, $h _{ij}=h(\breve{e} _i,\breve{e} _j)$
satisfy
\begin{equation}\label{decay condition}
\begin{aligned}
 a_{ij}=O(e^{- \tau \kappa r}),& \ \breve{{\nabla}}_k a_{ij}=O(e^{- \tau \kappa r}), \
\breve{\nabla}_l\breve{{\nabla}}_k a_{ij}=O(e^{-\tau \kappa r}),\\
h_{ij}=&O(e^{- \tau \kappa r}), \
\breve{{\nabla}}_k h_{ij}=O(e^{-\tau \kappa r}),
\end{aligned}
\end{equation}
where
$\breve{{\nabla}}$ is the Levi-Civita connection with respect to the hyperbolic metric
$\breve{g}$;\\
$(3)$  there is a distance function $\rho_z$
such that $T_{00} e^{\kappa \rho_z}$, $T_{0i} e^{\kappa \rho_z}$ $\in L^1(M)$. \end{defn}

\begin{rmk}
For simplicity, we assume the manifold $M$ has only one end. The results we obtain in the paper
could be extended to multi-end case easily.
\end{rmk}

Denote
$$\mathcal{E}_i=\breve{\nabla}^j g_{ij}-\breve {{\nabla}}_itr_{\breve{g}}(g)-\kappa(a_{1i}-g_{1i}tr_{\breve{g}}(a)),$$
$$\mathcal{P}_{ki}=p_{ki}-g_{ki}tr_{\breve{g}}(p).$$
Then we can define the following quantities for asymptotically AdS initial data.

\begin{defn}\label{definition}
For asymptotically $AdS$ initial data, the total energy is defined as
\begin{equation*}
E_0=\frac{\kappa}{16\pi}\lim_{r\rightarrow \infty}\int_{S_r}\mathcal{E}_1 U_{50}^{(0)}\breve{\omega}.
\end{equation*}
The total momenta are defined as
\begin{equation*}
\begin{aligned}
c_{i}=&\frac{\kappa}{16\pi}\lim_{r\rightarrow \infty}\int_{S_r}\mathcal{E}_1 U_{i5}^{(0)}\breve{\omega},\\
c'_{i}=&\frac{\kappa}{8\pi}\sum_{j=2}^{4}\lim_{r\rightarrow \infty}\int_{S_r}\mathcal{P}_{j1}U_{i0}^{(j)} \breve{\omega},\\
 J_{ij}=&\frac{\kappa}{8\pi}\sum_{j=2}^{4}\lim_{r\rightarrow \infty}\int_{S_r}\mathcal{P}_{k1}U_{ij}^{(k)} \breve{\omega},
\end{aligned}
\end{equation*}
where     $ \  \ \  \ \  \ \  \ \  \ \  \ \
\ \  \ \  \  \breve{\omega}= \breve{e}^2\wedge \breve{e}^3\wedge \breve{e}^4,\  \  \  \ U_{\alpha\beta}=U_{\alpha\beta}^{(\gamma)}\breve{e}_{\gamma}.$
\end{defn}

\begin{rmk}
If  $\kappa=1$, similar to \cite{WXZ}, we can derive the following relations between
the quantities in \cite{CMT} and the quantities defined in Definition $\ref{definition}$:
 \beQ
 \begin{aligned}
H(V_{(0)},0)&=m_{(0)}=E_0, \ \ \  \ H(V_{(i)},0)=m_{(i)}=-c_i,\\
H(0,C_{(i)})&=c_{(i)}=c_{i}',\ \ \ \ H(0,\Omega_{(i)(j)})=J_{(i)(j)}=J_{ij},
 \end{aligned}
 \eeQ
where  $i, j=1,2,3,4.$
\end{rmk}

\mysection{Positive energy theorem}\ls

Suppose $(N,\widetilde{g})$ is a $(4+1)$-dimensional spacetime, and $M$ is an asymptotically AdS hypersurface in $N$
with the induced metric $g$ and the second fundamental form $p$. $\widetilde{\nabla}$ and $\nabla$ are the Levi-Civita
connections corresponding to $\widetilde{g}$ and $g$ respectively. For simplicity, we also use the same symbols to
denote their lifts to the spinor bundle $\mathbb{S}$ respectively.
Define
\begin{equation}\label{Weitz formula}
\widehat{\nabla}_i=\widetilde{\nabla}_i+\frac{\sqrt{-1}}{2}\kappa
e_i\cdot,\ \ \widehat{D}=\sum_{i=1}^{4}e_i\cdot \widehat{\nabla}_i,
\end{equation}
then we can derive the following Weitzenb\"{o}ck formula \cite{XZ}
\begin{equation*}
\widehat{D}^*\widehat{D}=\widehat{\nabla}^{*}\widehat{\nabla}+\widehat{\mathcal{R}},
\end{equation*}
with
\begin{equation*}
\widehat{\mathcal{R}}=\frac{1}{2}(T_{00}e_0+T_{0i}e_i)\cdot e_0\cdot .
\end{equation*}
By Lax-Milgram Theorem, it is easy to prove that there exists a unique solution to the
equation $\widehat{D}\phi=0$ on $M$, with $\phi$ asymptotical to the imaginary Killing
spinor $\Phi_0$ on the end \cite{XZ}.
By integrating the Weitzenb\"{o}ck formula (\ref{Weitz formula}) and applying Witten's argument \cite{M,AnD,Z1,XZ}, we get
\begin{equation}\label{integral form of Weitzenbock formula}
\begin{aligned}
&\int_M|\widehat\nabla\phi|^2\ast 1
 +\int_M\langle\phi,\widehat{\mathcal{R}}\phi\rangle\ast1\\
=&\ \lim_{r\rightarrow
 \infty}Re\int_{S_r}\langle\phi,\sum_{ j\neq i }e_i\cdot
 e_j\cdot\widehat\nabla_j \phi\rangle\ast e^i\\
=&\ \frac{1}{4}\lim_{r\rightarrow
 \infty}\int_{S_r}(\breve{\nabla}^j
 g_{1j}-\breve
 {{\nabla}}_1tr_{\breve{g}}(g))|\Phi_0|^2\breve{\omega}\\
&\ +\frac{1}{4}\lim_{r\rightarrow
 \infty}\int_{S_r} \kappa(a_{k1}-g_{k1}tr_{\breve{g}}(a))\langle\Phi_0,\sqrt{-1}\breve{e}_k\cdot\Phi_0
 \rangle\breve{\omega}\\
&\ -\frac{1}{2}\lim_{r\rightarrow
 \infty}\int_{S_r}(h_{k1}-g_{k1}tr_{\breve{g}}(h))\langle\Phi_0,\breve{e}_0\cdot\breve{e}_k\cdot\Phi_0\rangle
 \breve{\omega},
\end{aligned}
\end{equation}
where $\phi$ is the unique solution of the equation $\widehat{D}\phi=0$.

By the Clifford representation (\ref{repre}) and the explicit form (\ref{ik}) of $\Phi_0$, the boundary term
on the right hand side of (\ref{integral form of Weitzenbock formula}) is equal to
\begin{equation}\label{RHS of integral form of Weitzenbock formula}
\begin{aligned}
RHS=&\
\frac{1}{2}\lim_{r\rightarrow\infty}\int_{S_r}\mathcal{E}_1 \Big(\overline{u^+}u^+ +\overline{v^+}v^+\Big) e^{\kappa r}\breve{\omega}\\
&\ +\lim_{r\rightarrow\infty}\int_{S_r}\mathcal{P}_{21}\Big(\overline{u^+}u^+ -\overline{v^+}v^+\Big)e^{\kappa r}\breve{\omega}\\
&\ -\sqrt{-1}\lim_{r\rightarrow\infty}\int_{S_r}\mathcal{P}_{31}\Big(\overline{u^+}v^+ -\overline{v^+}u^+\Big)e^{\kappa r}\breve{\omega}\\
&\ +\lim_{r\rightarrow\infty}\int_{S_r}\mathcal{P}_{41}\Big(\overline{u^+}v^+ +\overline{v^+}u^+\Big)e^{\kappa r}\breve{\omega}\\
=&\ 8\pi (\bar\lambda_1, \bar\lambda_2, \bar\lambda_3, \bar\lambda_4)Q(\lambda_1, \lambda_2, \lambda_3, \lambda_4)^t,
\end{aligned}
\end{equation}
in which the matrix
\begin{equation*}
Q=\begin{pmatrix}
E        &     L \\
\bar{L}^t  &    \hat{E}
\end{pmatrix},
\end{equation*}
where
\begin{equation*}
E=\begin{pmatrix}
E_0+c_4                   &      c_1'+\sqrt{-1}c_2'\\
+c_3'-J_{34}              &      -J_{14}-\sqrt{-1}J_{24}\\
\ & \ \\
c_1'-\sqrt{-1}c_2'        &      E_0+c_4\\
-J_{14}+\sqrt{-1}J_{24}   &      -c_3'+J_{34}
\end{pmatrix},
\end{equation*}

\begin{equation*}
\hat{E}=\begin{pmatrix}
E_0-c_4                   &      -c_1'-\sqrt{-1}c_2'\\
-c_3'-J_{34}              &      -J_{14}-\sqrt{-1}J_{24}\\
\ & \ \\
-c_1'+\sqrt{-1}c_2'       &      E_0-c_4\\
-J_{14}+\sqrt{-1}J_{24}   &      +c_3'+J_{34}
\end{pmatrix},
\end{equation*}

\begin{equation*}
L=\begin{pmatrix}
c_3                       &      c_1+\sqrt{-1}c_2\\
-c_4'+\sqrt{-1}J_{12}              &      +J_{13}+\sqrt{-1}J_{23}\\
\ & \ \\
c_1-\sqrt{-1}c_2          &      -c_3\\
-J_{13}+\sqrt{-1}J_{23}   &      -c_4'-\sqrt{-1}J_{12}
\end{pmatrix}.
\end{equation*}
Set $\hat{J}_k=\frac{1}{2}\varepsilon_{ijk}J_{jk}$ and denote
\beQ
\begin{aligned}
{\bf c}=&(c_1,c_2,c_3), \ {\bf c'}=(c'_1,c'_2,c'_3),\
{\bf\hat{J}}=(\hat{J}_1,\hat{J}_2,\hat{J}_3), \ {\bf J_{(4)}}=(J_{14},J_{24},J_{34}),\\
|L|^2=&2(|{\bf c}|^2+|{\bf\hat{J}}|^2+c_4'^2),
\ \ A=c_4^2+c_4'^2+|{\bf c}|^2+|{\bf c'}|^2+|{\bf\hat{J}}|^2+|{\bf J_{(4)}}|^2,
\end{aligned}
 \eeQ
then we have

\begin{thm} \label{main theorem}
Let $(M,g,h)$ be a $4$-dimensional asymptotically anti-de Sitter initial data of the spacetime
$(N,\widetilde{g})$ satisfying the dominant energy condition $(\ref{dec})$. Then we have the
following inequality:
\beq
E_0\geq max\Big \{&\big(c_4^2+\frac{1}{4}|L|^2\big)^\frac{1}{2}, \Big(\frac{1}{2}(|{\bf c}|^2+|{\bf J_{(4)}}|^2)+\frac{1}{8}|L|^2\Big)^\frac{1}{2},
(A+|{\bf c'}|^2+|{\bf J_{(4)}}|^2)^\frac{1}{2} \nonumber\\
&-|{\bf c'}|-|{\bf J_{(4)}}|,\ \ \big(A-2\sqrt{2}\big(\sum_{i=1}^3(c_4 c_{i}'-c_4' c_i)^2+|{\bf c} \times {\bf\hat{J}}|^2\big)^\frac{1}{2}\big)^\frac{1}{2},\nonumber\\
&\big(A-4\sqrt{2}\big(\sum_{i}(c_4 c_{i}'-c_4' c_i)^2+|{\bf c} \times {\bf\hat{J}}|^2\big)^\frac{1}{2}+F_+^\frac{1}{2}\big)^\frac{1}{2} \Big\},\nonumber
 \eeq
where
\beQ
 \begin{aligned}
F_+=max \{F,0\},
\end{aligned}
 \eeQ

\beQ
\begin{aligned}
F=&-8\sqrt{2}\big(\sum_{i=1}^3(c_4 c_{i}'-c_4' c_i)^2+|{\bf c} \times {\bf\hat{J}}|^2\big)^\frac{1}{2}A+36|{\bf c} \times {\bf\hat{J}}|^2+4|{\bf c} \times {\bf c'}|^2\\
&+36\sum_{i=1}^3(c_4 c_{i}'-c_4' c_i)^2
+4(|{\bf J_{(4)}}\cdot {\bf c'}|^2+ |{\bf J_{(4)}}\cdot {\bf \hat{J}}|^2+|{\bf J_{(4)}}\cdot {\bf c}|^2)\\
&+4|{\bf c'} \times {\bf \hat{J}}|^2+4|{\bf J_{(4)}}|^2(c_4^2+c_4'^2)+8c_4\varepsilon_{ijk} c_i \hat{J}_j J_{k4}
+8c_4'\varepsilon_{ijk}c_i' \hat{J}_{j} J_{k4}.
\end{aligned}
\eeQ
Moreover, if $E_0=0$, then $Q=0$ and the spacetime $(N,\widetilde{g})$ is anti-de Sitter along $M$. \end{thm}

\pf The nonnegativity of the Hermitian matrix $Q$ can be derived from the integral form of the Weitzenb\"{o}ck formula
 (\ref{integral form of Weitzenbock formula}), (\ref{RHS of integral form of Weitzenbock formula}) and the dominant energy condition $(\ref{dec})$.

The nonnegativity of first-order principal minors ensures $E_0\geq0$.
And from the nonnegativity of second-order principal minors, one finds
\begin{equation} \label{ine1}
E_0\geq \Big(c_4^2+\frac{1}{2}\sum_{i=1}^3\big(c_i^2+\hat{J}_i^2\big)+\frac{1}{2}c_4'^2\Big)^\frac{1}{2}
\end{equation}
 and
\begin{equation}\label{ine2}
E_0\geq \Big(\frac{1}{2}\sum_{i=1}^3\big(c_i'^2+J_{i4}^2\big)+\frac{1}{4}\sum_{i=1}^3\big(c_i^2+\hat{J}_i^2\big)+\frac{1}{4}c_4'^2\Big)^\frac{1}{2}.
\end{equation}

The sum of third-order principal minors is given, up to a positive constant, by
\beQ
 \begin{aligned}
S=
&E_0(E_0^2-A)+2c_4'\sum_{i=1}^3c_i J_{i4}+2\varepsilon_{ijk}c_i c_{j}' \hat{J}_k-2c_4\sum_{i=1}^3c_i'J_{i4}.
\end{aligned}
 \eeQ
Using the Cauchy inequality, one derives
\beQ
 \begin{aligned}
S\leq
&E_0(E_0^2-A)+2|c_4'||{\bf c}||{\bf J_{(4)}}|+2|{\bf c'}||{\bf c}\times {\bf \hat{J}}|+2E_0|{\bf c'}||{\bf J_{(4)}}|.
\end{aligned}
 \eeQ
Since $|c_4'||{\bf c}|\leq\frac{1}{2}(c_4'^2+|{\bf c}|^2)\leq E_0^2$
and
$|{\bf c}\times {\bf \hat{J}}|\leq\frac{1}{2}(|{\bf c}|^2+|{\bf \hat{J}}|^2)\leq E_0^2$,
one can obtain
\beQ
 \begin{aligned}
E_0\big((E_0+|{\bf c'}|+|{\bf J_{(4)}}|)^2-A-|{\bf c'}|^2-|{\bf J_{(4)}}|^2\big)\geq0.
\end{aligned}
 \eeQ
When $E_0>0$, we get
\beQ
\begin{aligned}
E_0\geq(A+|{\bf c'}|^2+|{\bf J_{(4)}}|^2)^\frac{1}{2}-|{\bf c'}|-|{\bf J_{(4)}}|,
\end{aligned}
 \eeQ
as
$A+|{\bf c'}|^2+|{\bf J_{(4)}}|^2\geq(|{\bf c'}|+|{\bf J_{(4)}}|)^2$.
When $E_0=0$, the inequality $(\ref{ine1})$, together with the inequality $(\ref{ine2})$, shows that $Q=0$. In this case, the inequality becomes trivial.

Also, we have
\beQ
\begin{aligned}
 &2c_4'\sum_{i=1}^3c_i J_{i4}+2\varepsilon_{ijk}c_i c_{j}' \hat{J}_k-2c_4\sum_{i=1}^3c_i'J_{i4}\\
 \leq &2(|{\bf J_{(4)}}|^2+|{\bf c'}|^2)^\frac{1}{2}\Big(\sum_{i=1}^3(c_4 c_{i}'-c_4' c_i)^2+|{\bf c} \times {\bf \hat{J}}|^2\Big)^\frac{1}{2}\\
\leq &2\sqrt{2}E_0\big(\sum_{i=1}^3(c_4 c_{i}'-c_4' c_i)^2+|{\bf c} \times {\bf \hat{J}}|^2\big)^\frac{1}{2}.
\end{aligned}
\eeQ
Therefore,
\beQ
 \begin{aligned}
S\leq
&E_0(E_0^2-A)+2\sqrt{2}E_0\big(\sum_{i=1}^3(c_4 c_{i}'-c_4' c_i)^2+|{\bf c} \times {\bf \hat{J}}|^2\big)^\frac{1}{2}.
\end{aligned}
 \eeQ
This implies
\beQ
\begin{aligned}
E_0^2\geq A-2\sqrt{2}\big(\sum_{i=1}^3(c_4 c_{i}'-c_4' c_i)^2+|{\bf c} \times {\bf \hat{J}}|^2\big)^\frac{1}{2}
\end{aligned}
\eeQ
if $E_0>0$. The case for $E_0=0$ is considered similarly.

The determinant of the matrix is
\beQ
\begin{aligned}
det Q=&\big(E_0^2-A\big)^2
+8E_0\sum_{i=1}^3(c_4'c_i J_{i4}-c_4c_i'J_{i4})+8E_0\varepsilon_{ijk}c_i c_{j}' \hat{J}_k\\
&-4|{\bf c} \times {\bf c'}|^2-4|{\bf c} \times {\bf \hat{J}}|^2-4|{\bf c'} \times {\bf \hat{J}}|^2-4|{\bf J_{(4)}}|^2(c_4^2+c_4'^2)\\
&-4\sum_{i=1}^3(c_4 c_{i}'-c_4' c_i )^2-4(|{\bf J_{(4)}}\cdot {\bf c'}|^2+ |{\bf J_{(4)}}\cdot {\bf \hat{J}}|^2+|{\bf J_{(4)}}\cdot {\bf c}|^2)\\
&-8c_4\varepsilon_{ijk} c_i \hat{J}_j J_{k4}
-8c_4'\varepsilon_{ijk}c_i' \hat{J}_{j} J_{k4}.
\end{aligned}
 \eeQ
Since
\beQ
\begin{aligned}
c_4'\sum_{i=1}^3c_i J_{i4}+\varepsilon_{ijk}c_i c_{j}' \hat{J}_k-c_4\sum_{i=1}^3c_i'J_{i4}
\leq \sqrt{2}E_0\big(\sum_{i=1}^3(c_4 c_{i}'-c_4' c_i)^2+|{\bf c} \times {\bf \hat{J}}|^2\big)^\frac{1}{2},
\end{aligned}
\eeQ
one obtains
\beQ
 \begin{aligned}
det Q
\leq&\big(E_0^2-A\big)^2
+8\sqrt{2}E_0^2\big(\sum_{i=1}^3(c_4 c_{i}'-c_4' c_i)^2+|{\bf c} \times {\bf \hat{J}}|^2\big)^\frac{1}{2}\\
&-4|{\bf c} \times {\bf c'}|^2-4|{\bf c} \times {\bf \hat{J}}|^2-4|{\bf c'} \times {\bf \hat{J}}|^2-4|{\bf J_{(4)}}|^2(c_4^2+c_4'^2)\\
&-4\sum_{i=1}^3(c_4 c_{i}'-c_4' c_i )^2-4(|{\bf J_{(4)}}\cdot {\bf c'}|^2+ |{\bf J_{(4)}}\cdot {\bf \hat{J}}|^2+|{\bf J_{(4)}}\cdot {\bf c}|^2)\\
&-8c_4\varepsilon_{ijk} c_i \hat{J}_j J_{k4}
-8c_4'\varepsilon_{ijk}c_i' \hat{J}_{j} J_{k4}\\
=&\Big(E_0^2-A+4\sqrt{2}\big(\sum_{i=1}^3(c_4 c_{i}'-c_4' c_i)^2+|{\bf c} \times {\bf \hat{J}}|^2\big)^\frac{1}{2}\Big)^2\\
&+8\sqrt{2}\big(\sum_{i=1}^3(c_4 c_{i}'-c_4' c_i)^2+|{\bf c} \times {\bf \hat{J}}|^2\big)^\frac{1}{2}A-36|{\bf c} \times {\bf \hat{J}}|^2-4|{\bf c} \times {\bf c'}|^2\\
&-36\sum_{i}(c_4 c_{i}'-c_4' c_i)^2-4(|{\bf J_{(4)}}\cdot {\bf c'}|^2+ |{\bf J_{(4)}}\cdot {\bf \hat{J}}|^2+|{\bf J_{(4)}}\cdot {\bf c}|^2)\\
&-4|{\bf c'} \times {\bf \hat{J}}|^2-4|{\bf J_{(4)}}|^2(c_4^2+c_4'^2)-8c_4\varepsilon_{ijk} c_i \hat{J}_j J_{k4}
-8c_4'\varepsilon_{ijk}c_i' \hat{J}_{j} J_{k4}.
\end{aligned}
 \eeQ
This implies
\beQ
E_0^2\geq A-4\sqrt{2}\big(\sum_{i}(c_4 c_{i}'-c_4' c_i)^2+|{\bf c} \times {\bf \hat{J}}|^2\big)^\frac{1}{2}+F_+^\frac{1}{2}.
\eeQ
The inequality claimed in the theorem follows immediately.

The rigidity part can be proved by following the argument in \cite{WXZ}. Here we skip the details.
\qed

\begin{rmk}
If $c_i=0, i=1,2,3,4$ and $c'_2=c'_4=J_{13}=J_{23}=J_{14}=J_{24}=0$ after suitable coordinate transformation,
the inequality for $(4+1)$-dimensional case in Theorem $2$ of \cite{CMT} can be derived from Theorem $\ref{main theorem}$.
\end{rmk}

\section{Appendix}\ls

\beQ \begin{aligned}
U_{50}=&\kappa^{-1} \frac{\partial}{\partial t},\\
U_{10}=&\kappa^{-1}\sin\theta \sin\psi \cos\varphi \frac{\partial}{\partial r}
+\coth(\kappa r)\Big(\cos\theta \sin\psi \cos\varphi \frac{\partial}{\partial \theta} \\
&+\frac{\cos\psi \cos\varphi}{\sin\theta} \frac{\partial}{\partial \psi}- \frac{\sin\varphi}{\sin\theta \sin\psi} \frac{\partial}{\partial \varphi}\Big),\\
U_{20}=&\kappa^{-1}\sin\theta \sin\psi \sin\varphi \frac{\partial}{\partial r} +\coth(\kappa r)\Big(\cos\theta \sin\psi \sin\varphi \frac{\partial}{\partial \theta} \\
&+ \frac{\cos\psi \sin\varphi}{\sin\theta} \frac{\partial}{\partial \psi}+ \frac{\cos\varphi}{\sin\theta \sin\psi} \frac{\partial}{\partial \varphi}\Big),\\
U_{30}=&\kappa^{-1}\sin\theta \cos\psi \frac{\partial}{\partial r} +\coth(\kappa
r)\Big(\cos\theta \cos\psi \frac{\partial}{\partial \theta}- \frac{\sin\psi}{\sin\theta} \frac{\partial}{\partial \psi}\Big),\\
U_{40}=&\kappa^{-1}\cos\theta \frac{\partial}{\partial r} -\coth(\kappa
r)\sin\theta\frac{\partial}{\partial \theta},\\
U_{15}=&\kappa^{-1}\tanh(\kappa r) \sin\theta \sin\psi \cos\varphi \frac{\partial}{\partial t},\\
U_{25}=&\kappa^{-1}\tanh(\kappa r) \sin\theta \sin\psi \sin\varphi \frac{\partial}{\partial t},\\
U_{35}=&\kappa^{-1}\tanh(\kappa r) \sin\theta \cos\psi \frac{\partial}{\partial t},\\
U_{45}=&\kappa^{-1}\tanh(\kappa r) \cos\theta \frac{\partial}{\partial t},\\
U_{12}=&\frac{\partial}{\partial \varphi},\\
U_{13}=&-\cos\varphi\frac{\partial}{\partial \psi}+\frac{\cos\psi \sin\varphi}{\sin \psi}\frac{\partial}{\partial \varphi},\\
U_{14}=&-\sin\psi \cos\varphi \frac{\partial}{\partial \theta}-\frac{\cos\theta \cos\psi \cos\varphi}{\sin \theta}\frac{\partial}{\partial \psi}+\frac{\cos\theta \sin\varphi}{\sin \theta \sin\psi}\frac{\partial}{\partial \varphi},\\
U_{23}=&-\sin\varphi\frac{\partial}{\partial \psi}-\frac{\cos\psi \cos \varphi}{\sin \psi}\frac{\partial}{\partial \varphi},\\
U_{24}=&-\sin\psi \sin\varphi \frac{\partial}{\partial \theta}-\frac{\cos\theta \cos\psi \sin\varphi}{\sin \theta}\frac{\partial}{\partial \psi}-\frac{\cos\theta \cos\varphi}{\sin \theta \sin\psi}\frac{\partial}{\partial \varphi},\\
U_{34}=&-\cos\psi \frac{\partial}{\partial \theta}+\frac{\cos\theta \sin\psi}{\sin \theta}\frac{\partial}{\partial \psi}.
\end{aligned} \eeQ

\end{document}